\documentclass{llncs}
\usepackage{amsmath,amssymb}
\usepackage{graphicx}

\begin{document}
\title{Quadri-allele frequency spectrum in a coalescent topology for mutations in non-constant population size}

\subtitle{}

\author{Arka Bhattacharya}

\institute{Columbia University,\ New York\ \email{\it (ab3899@columbia.edu)}}

\maketitle

\begin{abstract}
The sample frequency spectrum of a segregating site is the probability distribution of a sample of alleles from a genetic locus, conditional on observing the sample to have more than one clearly different phenotypes. We present a model for analyzing quadri-allele frequency spectrum, where the ancestral population diverged into three populations at a certain divergence time and the resulting mutations on the branches of the coalescent tree gave rise to three different derived alleles, which could be observed in the present generation along with the ancestral allele. The model has been analyzed for non-constant population size, assuming we had a certain number of extant lineages at the divergence time and no migration occurs between the populations.
\end{abstract}

\section{Introduction}
The sample frequency spectrum is basically the probability distribution of the number of mutant alleles in a sampled population, of the current generation.\ Since, the development of the elegant Coalescent Theory by Kingman[6], the analysis of the genealogical trees has become more strong and more mathematical in its approach, thus providing many information about allele frequencies, evolutionary histories etc.\ Assuming the infinite site model of mutation,\ the sample frequency spectrum is known in closed form, where the model is standard, neutral and has a coalescent approach\ (Griffiths and Tavare [7]).\ Recent works of Jenkins et. al [1],[2] obtain general trialleleic frequency spectrum using a coalescent approach, for varying population sizes.\ They obtained the predicted frequency spectrum of a site in a closed form, such that the site experienced at most two mutation events.\ In this work, we basically tried to obtain a quadri-allele frequency spectrum, using almost the same assumptions as that of [1] and [2], but the analysis and the construction of the model is different from them, though it incorporates one of their main results. \newline \newline 
The paper is divided into three sections.\ Section 2 gives a short overview of the work using pictorial representations.\ Section 3 gives the total mathematical analysis of the model in a very detailed manner and finally, Section 4 is the conclusion.

\section{Short explanation of the proposed work}
The work basically, tries to find the quadri-allele frequency spectrum of a non-constant population with a known most common recent ancestor (MRCA), using a coalescent approach, such that we see four different alleles in the present generation, including the ancestral one.\ Since, there must be at least three mutations in the genealogical tree to give rise to three different derived alleles, hence we considered a topology, where the ancestral population diverged into three different populations at a certain time and at least one non-nested mutation hit each of them.\ We consider the following general topology, given in Fig. 1.\ The green, blue and the yellow nodes represent the alleles, $b$, $c$ and $d$ respectively and each of the bold marked boxes, represents a diverged population, which diverged at a certain divergence time.\ To, have deeper understanding of the tree, consider a specific example of the topology in Fig. 2.\newline

\begin{figure}[ht!]
\centering
\includegraphics[width=125mm]{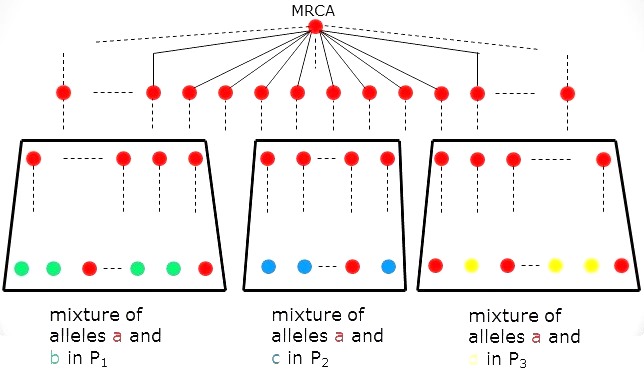}
\caption{General Topology}
\end{figure}

\begin{figure}[ht!]
\centering
\includegraphics[width=125mm]{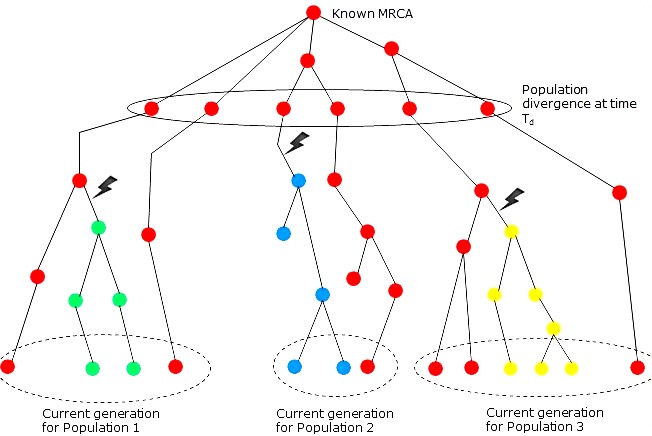}
\caption{Specific Topology}
\end{figure}
Now, the populations are combined in groups of two, in all the three possible ways and each such group is refered as a model. The expected time lengths for each of the models are found out and using that, the corresponding probability distribution for all of them are calculated and finally combined under an exponential family weight function, so that a model is incorporated in the final calculation, if and only if the allele frequency spectrum crosses a certain threshold.\ Refer, Fig. 3 for a pictorial representation of the division of the populations into separate models, for analysis.
\begin{figure}[ht!]
\centering
\includegraphics[width=70mm]{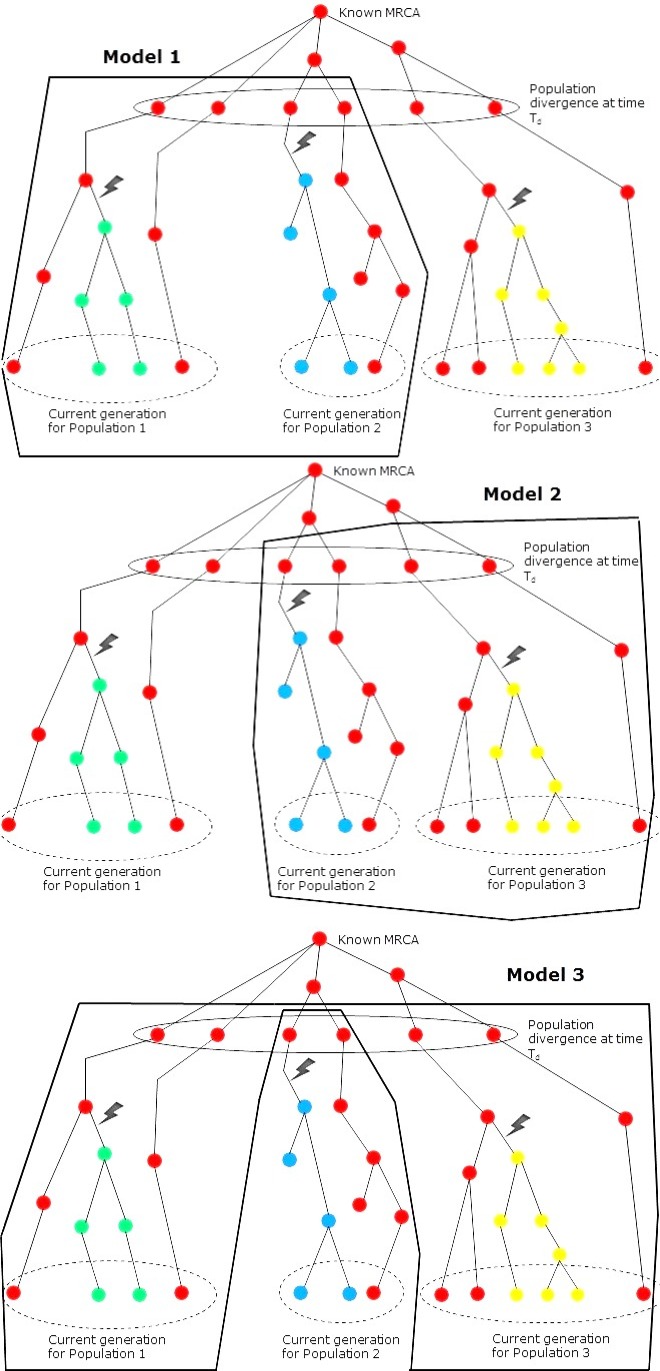}
\caption{Division into models}
\end{figure}
\section{Detailed mathematical analysis}
From, the previous diagrams, we saw that the ancestral population diverged into three different populations during the time of divergence and mutations occurred in each of those populations, only after divergence and gave rise to three different derived alleles in the present generation.\ In this section, we provide the overall mathematical details of the probability distribution of the allele frequency, conditioned on the event that we have sampled some number of individuals in the current generation, having a mixture of four different alleles. (one ancestral and three derived)\newline
Let us consider that $n$ lineages have been sampled in the current generation, consisting of $n_{1}$, $n_{2}$ and $n_{3}$ samples taken respectively from the first, second and the third population. We denote the respective populations by $P_{1}$, $P_{2}$ and $P_{3}$. Let the extant population present in all the three populations, during divergence be $m_{1}$, $m_{2}$ and $m_{3}$ correspondingly, such that the net extant population during divergence is $m=m_{1}+m_{2}+m_{3}$.All the three possible combinations of two populations are taken into account and are divided into separate models for analysis of allelular frequency. Model 1 represents $P_{1}$ and $P_{2}$ together, similarly model 2 and model 3 represent $P_{2},P_{3}$ and $P_{1},P_{3}$ together respectively.\ $M_{1},M_{2}$ and $M_{3}$ represent the extant lineages present in model 1,2 and 3 respectively, where we have $M_{1}=m_{1}+m_{2}$,\ $M_{2}=m_{2}+m_{3}$ and $M_{3}=m_{1}+m_{3}$.\ The number of sampled lineages in the current generation in model 1, 2 and 3 is $N_{1},N_{2}$ and $N_{3}$ respectively, where $N_{1}=n_{1}+n_{2}$, $N_{2}=n_{2}+n_{3}$ and $N_{3}=n_{1}+n_{3}$.\ If $n_{1}^{a}$ is the number of $a$ alleles found in $P_{1}$ and $n_{2}^{a}$ is the number of $a$ alleles found in $P_{2}$, then the number of $a$ alleles found in model 1 is $N_{1}^{a}=n_{1}^{a}+n_{2}^{a}$. Also, the number of $b$ and $c$ alleles in model 1, is $N_{1}^{b}=n_{1}^{b}$ and $N_{1}^{c}=n_{1}^{c}$. We must remember that $n_{1}=n_{1}^{a}+n_{1}^{b}$ and $n_{2}=n_{2}^{a}+n_{2}^{c}$ and for model 1, we have $N_{1}=N_{1}^{a}+N_{1}^{b}+N_{1}^{c}$.\ $T_{i}$ denotes the time in the coalescent tree at that point, where we have $i$ lineages present.\ Hence, $T_{N_{i}}$ and $T_{M_{i}}$ denote the starting and ending times, for the analysis of our models ($\forall i \in [1,3]$).\ Also, $P(t)$ denotes the number of lineages at time, $t$ for each of our models.\newline 
Thus, we have to find the expected time lengths of each stage for each of the models (\ finding for one suffices and then we can replace the corresponding parameters to find the values of the other two), find the corresponding probabilities of mutation and then have a frequency distribution of the alleles.\ Without loss of generality, we assume model 1 for our analysis.\ Hence, under model 1, we have to find the probability distribution of $N_{1}=N_{1}^{a}+N_{1}^{b}+N_{1}^{c}$, conditioned on the event that we have $M_{1}$ extant individuals present during divergence. (We assume that, physically the divergence time is equal for all the models and hence, to say that the number of lineages present at the divergence time, $t_{d}$ is $M_{1}$, we simply say $P(t_{d})=M_{1}$)\ Tavare [3] showed that, for constant population, the probability of the number of lineages at the time of divergence to be a certain value, can be found out. \newline 
Using his result, we have \newline 
\begin{equation*}
\hspace*{-9.5cm} Pr[P(t_{d})=M_{1}] = 
\end{equation*}
\begin{equation*}
\hspace*{-0.5cm} \sum_{i=M_{1}}^{N_{1}}\frac{(-1)^{i-M_{1}}(2i-1)[M_{1}.(M_{1}+1)...(M_{1}+i-2)].[N_{1}.(N_{1}-1)...(N_{1}-i+1)]}{M_{1}!(i-M_{1}).[N_{1}.(N_{1}+1)...(N_{1}+i-1)]}.e^{-\frac{i(i-1)}{4N}T_{M_{1}}}
\end{equation*}
\newline 
where, $N$ is the population size and represents the stationary allele frequency spectrum of equilibrium populations.\ The huge constant above, for each $i$, could be replaced by $c^{N_{1},M_{1},i}$.\ For general values of $N$ and $M$, we have \newline
\begin{equation*}
\hspace*{-11cm} c^{N,M,i} = 
\end{equation*}
\begin{equation*}
\hspace*{-0.2cm} \frac{(-1)^{i-M}(2i-1)[M.(M+1)...(M+i-2)].[N.(N-1)...(N-i+1)]}{M!(i-M).[N.(N+1)...(N+i-1)]}
\end{equation*}
\newline
so, that we have
\begin{equation*}
Pr[P(t_{d})=M_{1}]=\sum_{i=M_{1}}^{N_{1}}c^{N_{1},M_{1},i}.e^{-\frac{i(i-1)}{4N}T_{M_{1}}}
\end{equation*}
\newline 
Since, we are considering non-constant population size, we need to modify the above equation.\ Consider, the ratio of the present population size to the population at some time $t$, in a certain time frame, to be equal to $r(t)$.\ The above equation could be modified to,\newline 
\begin{equation*}
Pr[P(t_{d})=M_{1}]=\sum_{i=M_{1}}^{N_{1}}c^{N_{1},M_{1},i}.e^{-\frac{i(i-1)}{4N}\left(\int_{T_{N_{1}}}^{T_{N_{1}-1}}r(t)dt+\int_{T_{N_{1}-1}}^{T_{N_{1}-2}}r(t)dt+...+\int_{T_{M_{1}+1}}^{T_{M_{1}}}r(t)dt\right)}
\end{equation*}
\newline 
Now, since we consider, $T_{N_{1}}=0$, so the above equation starts from time 0.\ The above equation tells us the probability of the number of extant lineages at $t_{d}$ to be equal to $M_{1}$, when there are $N_{1}$ lineages in the present generation, so we could replace the quantity by the term $P^{N_{1},M_{1}}(t_{d})$.\ Hence, we have \newline
\begin{equation*}
P^{N_{1},M_{1}}(t_{d})=\sum_{i=M_{1}}^{N_{1}}c^{N_{1},M_{1},i}.\ e^{-\frac{i(i-1)}{4N}\int_{0}^{t_{d}}r(t)dt}
\end{equation*}
\newline
The above equation is also applicable for general values of $t$.\newline
Now, we need to find the expected values of time lengths,conditioned on the presence of $M_{1}$ extant lineages during divergence, i.e we need to evaluate $\mathbb{E}[T_{i}|P(t_{d})=M_{1}]$, for all values of $i \in [M_{1},N_{1}]$.\ Consider, the time taken for $i+1$ lineages to coalesce into $i$ lineages, be $T^{,}_{i}$ and let, $q^{T^{,}_{i}}(t)$ be the probability density function of $T^{,}_{i}$.\ Hence, if we need to find $\mathbb{E}[T_{i}|P(t_{d})=M_{1}]$, we must be able to find $q^{T^{,}_{i}|P(t_{d})=M_{1}}(t)$.\newline \newline 
For our first case, we will take $i$ to be equal to $M_{1}$ and hence, would obtain a formula for $\mathbb{E}[T_{M_{1}}|P(t_{d})=M_{1}]$.\ For infinitesimally small value of $h>0$, we have \newline 
\begin{equation*}
q^{T^{,}_{M_{1}}|P(t_{d})=M_{1}}(t)=\frac{1}{h}\lim_{h\rightarrow 0}Pr[t_{d}-t\leq T^{,}_{M_{1}}\leq t_{d}-t+h|P(t_{d})=M_{1}]
\end{equation*}
\begin{equation*}
=\frac{1}{P^{N_{1},M_{1}}(t_{d})}.\frac{1}{h}\lim_{h\rightarrow 0}Pr[P(t_{d}-t)=M_{1},P(t_{d}-t+h)=M_{1},P(t_{d})=M_{1}]
\end{equation*}
\begin{equation*}
\hspace*{-3.1cm}=\frac{1}{P^{N_{1},M_{1}}(t_{d})}.\frac{1}{h}\lim_{h\rightarrow 0}Pr[P(t_{d}-t)=M_{1},P(t_{d})=M_{1}]
\end{equation*}
\begin{equation*}
\hspace*{-2.3cm}=\frac{1}{P^{N_{1},M_{1}}(t_{d})}.P^{N_{1},M_{1}+1}(t_{d}-t).\frac{M_{1}.(M_{1}+1)}{4N}.e^{-\frac{M_{1}(M_{1}-1)}{4N}t}
\end{equation*}
\newline 
Expanding the above equation in terms of the known parameters, we finally get, 
\begin{equation*}
q^{T^{,}_{M_{1}}|P(t_{d})=M_{1}}(t)=\frac{\sum_{i=M_{1}+1}^{N_{1}}c^{N_{1},M_{1}+1,i}.\ e^{-\frac{i(i-1)}{4N}\int_{0}^{t_{d}-t}r(t)dt}}{\sum_{i=M_{1}}^{N_{1}}c^{N_{1},M_{1},i}.\ e^{-\frac{i(i-1)}{4N}\int_{0}^{t_{d}}r(t)dt}}.\frac{M_{1}.(M_{1}+1)}{4N}.e^{-\frac{M_{1}(M_{1}-1)}{4N}t}
\end{equation*}
Hence, we finally get,
\begin{equation*}
\mathbb{E}[T_{M_{1}}|P(t_{d})=M_{1}]=\int_{0}^{t_{d}}t.q^{T^{,}_{M_{1}}|P(t_{d})=M_{1}}(t) dt
\end{equation*}
\begin{equation*}
=\int_{0}^{t_{d}}t.\left(\frac{\sum_{i=M_{1}+1}^{N_{1}}c^{N_{1},M_{1}+1,i}.\ e^{-\frac{i(i-1)}{4N}\int_{0}^{t_{d}-t}r(t)dt}}{\sum_{i=M_{1}}^{N_{1}}c^{N_{1},M_{1},i}.\ e^{-\frac{i(i-1)}{4N}\int_{0}^{t_{d}}r(t)dt}}.\frac{M_{1}.(M_{1}+1)}{4N}.e^{-\frac{M_{1}(M_{1}-1)}{4N}t}\right)dt
\end{equation*}
\newline 
For our second case, we will consider, $M_{i}<i<N_{1}$.\newline 
Now, we have to consider all the possibilities of starting and ending points within the coalescent tree.\ To find the probability density distribution of $T^{,}_{i}$, we have to make sure that the population at the time of divergence is $M_{1}$, i.e $P(t_{d})=M_{1}$.\ Hence, if the starting point is at some time, $t_{1}$ and the ending point is at some time, $t_{d}-t$. then we must include the conditional probability of the event of seeing $M_{1}$ extant lineages during the divergence time, conditioned on the event that the coalescent time of $i-2$ lineages coalescing to $i-1$ lineages, is $t_{1}+t$.\ Thus, we have \newline 
\begin{equation*}
q^{T^{,}_{i}|P(t_{d})=M_{1}}(t)=\frac{1}{Pr[P(t_{d})=M_{1}]}.\left(\int_{t_{1}=0}^{t_{d}-t}q^{T^{,}_{i}}(t_{1}).q^{T_{i}}(t).Pr[P(t_{d})=M_{1}|T^{,}_{i-1}=t_{1}+t]dt_{1}\right)
\end{equation*}
\begin{equation*}
\hspace*{-0.5cm}=\frac{1}{P^{N_{1},M_{1}}(t_{d})}.\left(\int_{t_{1}=0}^{t_{d}-t}q^{T^{,}_{i}}(t_{1}).q^{T_{i}}(t).Pr[P(t_{d})=M_{1}|T^{,}_{i-1}=t_{1}+t]dt_{1}\right)
\end{equation*}
\begin{equation*}
\hspace*{-0.5cm}=\frac{1}{\sum_{i=M_{1}}^{N_{1}}c^{N_{1},M_{1},i}.e^{-\frac{i(i-1)}{4N}\int_{0}^{t_{d}}r(t)dt}}.\left(\int_{t_{1}=0}^{t_{d}-t}q^{T^{,}_{i}}(t_{1}).q^{T_{i}}(t).Pr[P(t_{d})=M_{1}|T^{,}_{i-1}=t_{1}+t]dt_{1}\right)
\end{equation*}
\newline
Finally, we get,\newline 
\begin{equation*}
\hspace*{-9.5cm}\mathbb{E}[T_{i}|P(t_{d})=M_{1}]=
\end{equation*}
\begin{equation*}
\int_{0}^{t_{d}}\left(\frac{t}{\sum_{i=M_{1}}^{N_{1}}c^{N_{1},M_{1},i}.e^{-\frac{i(i-1)}{4N}\int_{0}^{t_{d}}r(t)dt}}.\int_{t_{1}=0}^{t_{d}-t}q^{T^{,}_{i}}(t_{1}).q^{T_{i}}(t).Pr[P(t_{d})=M_{1}|T^{,}_{i-1}=t_{1}+t]dt_{1}\right)dt
\end{equation*}
\newline
Now, we have our final case, when $i=N_{1}$. This case is similar to the case of $i=M_{1}$, instead
,we will have terms, in the parameter of $N_{1}-1$, instead of $M_{1}+1$ as in the first case.\ As usual, for infinitesimally small $h>0$, we have \newline 
\begin{equation*}
q^{T^{,}_{N_{1}}|P(t_{d})=M_{1}}(t)=\frac{1}{h}\lim_{h\rightarrow 0}Pr[t_{d}-t\leq T^{,}_{N_{1}-1}\leq t_{d}-t+h|P(t_{d})=M_{1}]
\end{equation*}
\begin{equation*}
\hspace*{-3cm}=\frac{1}{P^{N_{1},M_{1}}(t_{d})}.\ q^{T^{,}_{N_{1}-1}}(t).Pr[P(t_{d})=M_{1}|T^{,}_{N_{1}-1}=t]
\end{equation*}
\begin{equation*}
\hspace*{-2.2cm}=\frac{1}{P^{N_{1},M_{1}}(t_{d})}.{P^{N_{1}-1,M_{1}}(t_{d}-t).\frac{N_{1}(N_{1}-1)}{4N}.e^{-\frac{N_{1}(N_{1}-1)}{4N}t}}
\end{equation*}

Expanding, the above equation in terms of known parameters, we finally obtain \newline 
\begin{equation*}
q^{T^{,}_{N_{1}}|P(t_{d})=M_{1}}(t)=\frac{\sum_{i=M_{1}}^{N_{1}-1}c^{N_{1}-1,M_{1},i}.\ e^{-\frac{i(i-1)}{4N}\int_{0}^{t_{d}-t}r(t)dt}}{\sum_{i=M_{1}}^{N_{1}}c^{N_{1},M_{1},i}.\ e^{-\frac{i(i-1)}{4N}\int_{0}^{t_{d}}r(t)dt}}.\frac{N_{1}(N_{1}-1)}{4N}.e^{-\frac{-N_{1}(N_{1}-1)}{4N}t}
\end{equation*}
\newline 
Finally, \newline 
\begin{equation*}
\mathbb{E}[T_{N_{1}}|P(t_{d})=M_{1}]=\int_{0}^{t_{d}}t.q^{T^{,}_{N_{1}}|P(t_{d})=M_{1}}(t) dt
\end{equation*}
\begin{equation*}
=\int_{0}^{t_{d}}t.\left(\frac{\sum_{i=M_{1}}^{N_{1}-1}c^{N_{1}-1,M_{1},i}.\ e^{-\frac{i(i-1)}{4N}\int_{0}^{t_{d}-t}r(t)dt}}{\sum_{i=M_{1}}^{N_{1}}c^{N_{1},M_{1},i}.\ e^{-\frac{i(i-1)}{4N}\int_{0}^{t_{d}}r(t)dt}}.\frac{N_{1}(N_{1}-1)}{4N}.e^{\frac{-N_{1}(N_{1}-1)}{4N}t}\right)dt
\end{equation*}
Hence, we found out the equations, for all the expected values of time lengths, for model 1, considering that the extant population was $M_{1}$ and the current generation lineage is $N_{1}$.\ Similarly, by just manipulating the parameters, we could find the expected time lengths for both model 2 and model 3, the equations remaining the same.\ Thus, we finally obtained, \newline 
\begin{equation*}
\mathbb{E}[T_{i}|P(t_{d})=M_{1}]\ =\ \mathbb{E}[T_{i}]_{model 1}\ \ \forall \ M_{1}\leq i\leq N_{1}
\end{equation*}
\begin{equation*}
\mathbb{E}[T_{i}|P(t_{d})=M_{2}]\ =\ \mathbb{E}[T_{i}]_{model 2}\ \ \forall \ M_{2}\leq i\leq N_{2}
\end{equation*}
\begin{equation*}
\mathbb{E}[T_{i}|P(t_{d})=M_{3}]\ =\ \mathbb{E}[T_{i}]_{model 3}\ \ \forall \ M_{3}\leq i\leq N_{3}
\end{equation*}
Now, we use the tri-allele frequency spectrum results from [1] and [2] for each of our models. \newline \newline 
Since, we already mentioned that for each of our populations, we need to have at least one segregating site, to give rise to four different alleles in the present generation, including the ancestral one, so we consider the non-nested model of their analysis, so that we have high chances of observing meaningful mutations in each of the two populations, for every model that we considered.\ Without loss of generality, we again take model 1, to apply the results of [1] and [2].\ Recall that, we represented the number of sampled individuals having $a$, $b$ and $c$ alleles, by $N_{1}^{a}$,\ $N_{1}^{b}$ and $N_{1}^{c}$ respectively, so that we get $N_{1}=N_{1}^{a}+N_{1}^{b}+N_{1}^{c}$.\ We also had $M_{1}$ extant individuals, having the $a$ alleles.\ So, we are considering the case for the event of non-nested mutations, which gave rise to the alleles, $b$ and $c$ and let us denote this event, by $E_{bc}$.\ Let $P$ be the mutation transition matrix, where for every row $i$ and column $j$, $P_{ij}$ denotes the probability of transition from allele $i$ to allele $j$.\ We take $\theta$ to be the population scaled mutation rate and $u$ to be the probability of a mutation event at a given locus per meiosis, so that $\theta =4N_{0}u$, where $N_{0}$ is the diploid effective population size.\ Then, finally applying the results, we get the following probability distribution for model 1, so that we have two different derived alleles, $b$ and $c$, for $\theta$ tending to 0. \newline \newline 
\begin{equation*}
\hspace*{-7cm}Pr[N_{1}=N_{1}^{a}+N_{1}^{b}+N_{1}^{c},\ E_{bc}]=
\end{equation*}
\begin{equation*}
\left[\frac{\theta^{2}}{4}P_{ab}P_{ac}\sum_{k=M_{1}+1}^{N_{1}^{a}+N_{1}^{b}+1}\sum_{j=M_{1}}^{k}\left(\sum_{l=j-2}^{k-2}\frac{\binom{N_{1}^{a}-1}{l-1}\binom{N_{1}^{b}-1}{k-l-2}\binom{k-j}{k-l-2}}{\binom{N_{1}-1}{k-1}\binom{k-1}{l+1}}\frac{j(j-1)}{1+\delta_{j,k}} \right) \mathbb{E}[T_{j}T_{k}]\right]+O(\theta^{3})
\end{equation*}
where $\delta_{j,k}$ is the Kronecker delta. \newline \newline
The probability of the event, $E_{bc}$ could also be found out by the equation below. \newline \newline 
\begin{equation*}
Pr[E_{bc}]=\left[ \frac{\theta^{2}}{4}P_{ab}P_{ac} \sum_{k=M_{1}+1}^{N_{1}}\sum_{j=M_{1}}^{k}\left(\left( k(j-1)-\frac{2\delta_{j,2}}{k-1}\right)\frac{1}{1+\delta_{j,k}} \right)\mathbb{E}[T_{j}T_{k}] \right]+O(\theta^{3})
\end{equation*}
Similarly, we get the corresponding probabilities, for the other two models.\ Though the formula remains constant for the other two models, but still we write the equations below, by only changing the parameters.
\newline \newline 
\begin{equation*}
\hspace*{-7cm}Pr[N_{2}=N_{2}^{a}+N_{2}^{c}+N_{2}^{d},\ E_{cd}]=
\end{equation*}
\begin{equation*}
\left[\frac{\theta^{2}}{4}P_{ac}P_{ad}\sum_{k=M_{2}+1}^{N_{2}^{a}+N_{2}^{c}+1}\sum_{j=M_{2}}^{k}\left(\sum_{l=j-2}^{k-2}\frac{\binom{N_{2}^{a}-1}{l-1}\binom{N_{2}^{c}-1}{k-l-2}\binom{k-j}{k-l-2}}{\binom{N_{2}-1}{k-1}\binom{k-1}{l+1}}\frac{j(j-1)}{1+\delta_{j,k}} \right) \mathbb{E}[T_{j}T_{k}]\right]+O(\theta^{3})
\end{equation*}
\newline
\begin{equation*}
Pr[E_{cd}]=\left[ \frac{\theta^{2}}{4}P_{ac}P_{ad} \sum_{k=M_{2}+1}^{N_{2}}\sum_{j=M_{2}}^{k}\left(\left( k(j-1)-\frac{2\delta_{j,2}}{k-1}\right)\frac{1}{1+\delta_{j,k}} \right)\mathbb{E}[T_{j}T_{k}] \right]+O(\theta^{3})
\end{equation*}
\newline \newline \newline 
\begin{equation*}
\hspace*{-7cm}Pr[N_{3}=N_{3}^{a}+N_{3}^{d}+N_{2}^{b},\ E_{db}]=
\end{equation*}
\begin{equation*}
\left[\frac{\theta^{2}}{4}P_{ad}P_{ab}\sum_{k=M_{3}+1}^{N_{3}^{a}+N_{3}^{d}+1}\sum_{j=M_{3}}^{k}\left(\sum_{l=j-2}^{k-2}\frac{\binom{N_{3}^{a}-1}{l-1}\binom{N_{3}^{d}-1}{k-l-2}\binom{k-j}{k-l-2}}{\binom{N_{3}-1}{k-1}\binom{k-1}{l+1}}\frac{j(j-1)}{1+\delta_{j,k}} \right) \mathbb{E}[T_{j}T_{k}]\right]+O(\theta^{3})
\end{equation*}
\newline
\begin{equation*}
Pr[E_{db}]=\left[ \frac{\theta^{2}}{4}P_{ad}P_{ab} \sum_{k=M_{3}+1}^{N_{3}}\sum_{j=M_{3}}^{k}\left(\left( k(j-1)-\frac{2\delta_{j,2}}{k-1}\right)\frac{1}{1+\delta_{j,k}} \right)\mathbb{E}[T_{j}T_{k}] \right]+O(\theta^{3})
\end{equation*}
\newline 
Thus, we finally got the expected time lengths for all the models, conditioned on the event of an extant population during divergence and using those results, could finally calculate the probability distribution of the sampled individuals in each of the models and also find out the probability of the mutation events for all of them.\newline 
Now, we define the following simple exponential family weight function, for each of the models, such that the weight of the model in giving a final contribution to the allele frequency spectrum is higher if and only if the model has higher probability of a mutation event, giving rise to two different alleles.\ This weighted analysis, could be called \textit{\textbf{The Weighted Coalescent Approach}}.\ We have, \newline \newline 
\begin{equation*}
Weight_{model 1}=w_{1}=\frac{e^{Pr[E_{bc}]}}{e^{Pr[E_{bc}]}+e^{Pr[E_{cd}]}+e^{Pr[E_{db}]}}
\end{equation*}
\begin{equation*}
Weight_{model 2}=w_{2}=\frac{e^{Pr[E_{cd}]}}{e^{Pr[E_{bc}]}+e^{Pr[E_{cd}]}+e^{Pr[E_{db}]}}
\end{equation*}
\begin{equation*}
Weight_{model 1}=w_{3}=\frac{e^{Pr[E_{db}]}}{e^{Pr[E_{bc}]}+e^{Pr[E_{cd}]}+e^{Pr[E_{db}]}}
\end{equation*}
\newline 
Again, it is possible that we could get a finite value of the probability distribution for each of the models, yet a very low value for a meaningful contribution to the net allele frequency spectrum.\ So, it is logical to assign indicator variables for the models, to check if the allele frequency spectrum crosses a certain small threshold, $\epsilon_{i}$, $\forall i\in [3]$. such that $0<\epsilon_{i}<1$.\ Thus, we have, \newline \newline 
\[I_{1} = \left\{
  \begin{array}{lr}
    1 & :Pr[N_{1},E_{bc}]\geq \epsilon_{1}\\
    0 & :Pr[N_{1},E_{bc}]< \epsilon_{1}
  \end{array}
\right.
\]
Similarly, we have \newline 
\[I_{2} = \left\{
  \begin{array}{lr}
    1 & :Pr[N_{2},E_{cd}]\geq \epsilon_{2}\\
    0 & :Pr[N_{2},E_{cd}]< \epsilon_{2}
  \end{array}
\right.
\]
\[I_{3} = \left\{
  \begin{array}{lr}
    1 & :Pr[N_{3},E_{db}]\geq \epsilon_{3}\\
    0 & :Pr[N_{3},E_{db}]< \epsilon_{3}
  \end{array}
\right.
\]
\newline
Now, that we got the weights and threshold indicators for each of the models, hence the final allele frequency spectrum could be written as below. \newline 
\begin{equation*}
Pr[n=n^{a}+n^{b}+n^{c}+n^{d}]=I_{1}.w_{1}.Pr[N_{1},E_{bc}]+I_{2}.w_{2}.Pr[N_{2},E_{cd}]+I_{3}.w_{3}.Pr[N_{3},E_{db}]
\end{equation*}
\section{Conclusion}
Though, we saw the above model for finding out the quadri allele frequency spectrum, in which the present sampled individuals have four different alleles, but I personally feel, that allele frequency spectrum analysis is far from being perfect, since it relies on a number of assumptions, that are in fact not enforced by mother nature herself.\ But, proper statistical methods, need to be employed to include many factors like migration, natural selection etc.\ This simple weighted coalescent approach may be extended and made more general to incorporate far more complex factors in estimating allele frequency spectrums.\ It would be interesting to see, if we could divide the coalescent tree, like this, for $k-1$ populations, into $\binom{k-1}{2}$ models, so that we get a general $k$-allele frequency spectrum and develop far more complex weight functions, for each of the models, to incorporate migration and natural selection.\ Or, we could divide the tree into far lesser models, if possible, such that each model has 2 or more than 2 populations combined within it, and then recurse on that model, which has more than 2 diverged populations, until we get simpler models and then finally combine them.\ Well, it is true, that mother nature does not play dice, but she seems to play hide and seek with us, hiding all the important questions and their solutions.

\bibliographystyle{splncs}

\end{document}